\documentclass[12pt,twoside]{amsart}
\usepackage{amssymb}

\nonstopmode

\numberwithin{equation}{section}
\newtheorem{theorem}{Theorem}[section]
\newtheorem{lemma}[theorem]{Lemma}
\newtheorem{proposition}[theorem]{Proposition}
\newtheorem{corollary}[theorem]{Corollary}

\theoremstyle{definition}
\newtheorem{definition}[theorem]{Definition} 
\newtheorem{remark}[theorem]{Remark}
\newtheorem{example}[theorem]{Example}

\newcommand\Ann{\operatorname{Ann}}
\newcommand\Supp{\operatorname{Supp}}
\newcommand\Ass{\operatorname{Ass}}
\newcommand\Hom{\operatorname{Hom}}
\newcommand\Ext{\operatorname{Ext}}

\newcommand\depth{\operatorname{depth}}

\newcommand{\HH}{H_{\mathfrak m}}

\newcommand{\ER}{\operatorname{Ext}_R}
\newcommand{\Proj}{\operatorname{Proj}}

\newcommand{\s}{\; | \;}
\newcommand{\mif}{\mbox{if} ~}

\newcommand{\ffi}{\varphi}

\newcommand{\cE}{{\mathcal E}}
\newcommand{\cF}{{\mathcal F}}
\newcommand{\cG}{{\mathcal G}}

\newcommand{\fm}{{\mathfrak m}}
\newcommand{\fa}{{\mathfrak a}}
\newcommand{\fc}{{\mathfrak c}}
\newcommand{\fp}{{\mathfrak p}}

\newcommand {\ZZ}{\mathbb{Z}}

\newcommand {\PP}{\mathbb{P}}

\newcommand{\tk}{\otimes_K}
\newcommand{\tr}{\otimes_R}

\begin{document}
\title[Bezout's theorem and Cohen-Macaulay modules]{Bezout's theorem and Cohen-Macaulay modules}

\author[J.\ Migliore, U.\ Nagel, C.\ Peterson]{J.\ Migliore, U.\
  Nagel, C.\ Peterson}
\address{Department of Mathematics,
        University of Notre Dame,
        Notre Dame, IN 46556,
        USA}
\email{Juan.C.Migliore.1@nd.edu}
\address{Fachbereich Mathematik und Informatik, Universit\"at-Gesamthochschule
Paderborn, D--33095 Paderborn, Germany}
\email{uwen@uni-paderborn.de}
\address{Department of Mathematics,
Washington University,
St. Louis, MO 63130-4899}
\email{peterson@math.wustl.edu}



\begin{abstract} We define very proper intersections of modules and
  projective subschemes. It turns out that  equidimensional locally
  Cohen-Macaulay modules intersect 
  very properly if and only if they intersect properly. We prove a Bezout
  theorem for modules which meet very properly. Furthermore, 
  we show for equidimensional sub\-schemes $X$ and $Y$: If they intersect
  properly in an arithmetically Cohen-Macaulay sub\-scheme of positive
  dimension   then $X$ and $Y$ are arithmetically Cohen-Macaulay. The
  module version 
  of this result implies splitting criteria for reflexive sheaves. 
\end{abstract}


\maketitle

\tableofcontents


 \section{Introduction} \label{intro}

In this note we present a unified approach in order to study the
intersection of
projective subschemes and tensor products of reflexive sheaves on projective
space over an arbitrary field. A crucial role is played by the new concept of a
very proper intersection.

There are two starting points for our considerations.  The first is the
following result of Huneke and Ulrich.

\begin{proposition} \label{intro_HU-lifting}
Let $X \subset \PP^n$ denote an equidimensional subscheme having dimension
at least $2$. Let $H$ be a hyperplane which does not contain any component
of $X$. Suppose that  $X \cap H$ is arithmetically Cohen-Macaulay. Then
$X$ is arithmetically Cohen-Macaulay. 
\end{proposition}

We want to generalize this statement to a result on the intersection of two
projective subschemes. Thus we consider the following problem.
\smallskip

\noindent
{\bf Problem:} Let $X, Y \subset \PP^n$ be equidimensional subschemes such
that  $\dim X \cap Y \geq 1$. Assume that 
$X \cap Y$ is arithmetically Cohen-Macaulay.
Which ``extra conditions'' imply that $X$ and $Y$ are arithmetically
Cohen-Macaulay?  
\smallskip

We assume that the dimension of $X \cap Y$ is at least one because $X \cap
Y$ is
automatically arithmetically Cohen-Macaulay if it is zero-dimensional. Still,
a priori it is not clear at all if the problem above has a reasonable answer.
Suppose it has. Then, a comparison with the result of Huneke and Ulrich
suggests
to look for some genericity condition which ensures that the subschemes $X$ and
$Y$ are in sufficiently general position to each other. Thus, we want to
assume
that $X$ and $Y$ meet properly. Recall that there is always the following
inequality 
$$
\dim X \cap Y \geq \dim X + \dim Y - n.
$$
If equality holds true then it is said that $X$ and $Y$ intersect properly.
However, the following simple example shows that ``proper intersection'' is
not the answer to our problem.

\begin{example} \label{intro_example} Let $T_1, T_2, T_3 \subset \PP^5$
be the three-folds defined by $(x_0, x_1), \; (x_2, x_3)$ and $(x_1 + x_2,
x_0 + x_3)$, respectively. Put $X = T_1 \cup T_2$ and $Y = T_3$. Then
$I_X + I_Y$ is a saturated ideal and $X$ and $Y$ intersect properly in an
arithmetically Cohen-Macaulay curve. However, $X$ is not arithmetically
Cohen-Macaulay.
\end{example}

The analysis of this example lead us to the concept of a ``very
proper intersection'' (cf.\ Definition \ref{very_proper_inter-def}).
Intuitively, two subschemes $X$ and $Y$ intersect very properly if not only
$X$ and $Y$ intersect properly but also all their cohomology modules do. We
show that the two conditions proper intersection and very proper
intersection  are equivalent for $X$ and $Y$  if they are both
equidimensional and locally Cohen-Macaulay.
As one of our goals, we prove that the condition very proper intersection is
an answer to our problem above.  

\begin{theorem} \label{intro-thm-CM-lifting} Let $X, Y \subset \PP^n$
  denote equidimensional subschemes which meet very properly in an
  arithmetically Cohen-Macaulay subscheme of positive dimension. Then
  $X$ and $Y$ are arithmetically Cohen-Macaulay. 
\end{theorem} 

It is interesting to see that the example above does not contradict this
result.  Observe that $X$ is not locally Cohen-Macaulay 
along the line $T_1 \cap T_2$ which is contained in $Y$. Thus $X$ and $Y$
do not meet very properly.

A natural idea for proving this result is to use Serre's
diagonal trick. This amounts to viewing the intersection of $X$ and $Y$ as
consecutive hyperplane sections of their join. Then, one is tempted to apply
 Proposition \ref{intro_HU-lifting}. Note, however, that this statement
 assumes  that the scheme considered is  equidimensional. 
 Furthermore,  observe that the occurring hyperplanes are not general because
 they have to contain the diagonal. In fact, the example above shows that 
 things can go wrong. 
In order to control these problems we study the
 cohomology of the subschemes involved. Our methods provide the following
 version of Bezout's theorem. 

\begin{theorem} \label{intro-Bezout} 
Let $X, Y \subset \PP^n$
  denote equidimensional subschemes which intersect very properly in a
  non-empty scheme. Then the intersection is equidimensional, and we have 
$$
\deg X \cap Y = \deg X \cdot \deg Y. 
$$
\end{theorem}

It is well-known that the degree relation is not true for arbitrary proper
intersections. Very proper intersections have better properties. The
collection of irreducible subschemes of the intersection, which contribute
to the Bezout 
number $\deg X \cdot \deg Y$, corresponds  exactly to the whole set of irreducible
components of the homogeneous ideal of $X \cap Y$. In
particular, $X \cap Y$ has no embedded components.

In the special case of a  hypersurface section,  the last result gives a
condition when the hypersurface section 
of an equidimensional subscheme remains equidimensional. Thus we may also
view Theorem \ref{intro-Bezout}  as a version of Bertini's theorem.  

Our second starting point for this note were results of Huneke and Wiegand
in \cite{Huneke-Wiegand-Math.Ann} and \cite{Huneke-Wiegand-Math.Scand} on tensor products of maximal modules. Here we think of  tensor products
as  intersections of modules. 
It turns out that the concept of a
very proper intersection can be extended to the case of graded modules and
sheaves 
on projective space. In fact,
the results mentioned above are special cases of statements about  tensor
products of modules.  Note also that two reflexive sheaves on projective
space always meet  properly. Their tensor product is not necessarily
reflexive. However, we will show that it is indeed reflexive provided the
sheaves intersect very properly. 

The paper is organized as follows. In Section \ref{sec-unmixedness}, we
introduce the notion of a slight modification. Then we use cohomological
methods in order to relate the unmixedness of a module $M$ to the
unmixedness of $M/I M$ where $I$ is a parameter ideal for $M$. 
Section \ref{Kunneth} is entirely devoted to compute the local cohomology
of the  join of two modules. The resulting  K\"unneth formulas are
essential for the rest of the paper. Very proper intersections are
introduced in Section \ref{Bezout}. There we also prove our version of
Bezout's theorem for modules. In Section \ref{lifting},  we combine all
these 
techniques in order 
to relate the Cohen-Macaulay property of the tensor product to its
factors. The consequences for reflexive sheaves are drawn in the final
section.  

Throughout this note,   $R$ denotes a  standard graded $K$-algebra or a local
ring containing a field $K$. 
If $M$ is a module over the graded ring $R$ then it is understood that $M$
is a graded $R$-module
\medskip

\noindent
{\bf Acknowledgement.} The authors would like to thank Craig Huneke for
pointing out Serre's use of the diagonal trick.

Most of the present work was accomplished  during a stay of the second author
at the University of Notre Dame. He would like to thank its Department of
Mathematics for support and hospitality.


\section{Unmixedness results} \label{sec-unmixedness} 

In this section we assume that $R$ is a Gorenstein ring which is either
local or a graded $K$-algebra where $K$ is an infinite field.
In \cite{N-gorliaison}  the unmixedness of an ideal of $R$ has
been characterized by means of its local cohomology modules. First we note
that this result can be extended to finitely generated $R$-modules. Second
we apply it in order to derive sufficient conditions on a subsystem
$\{l_1,\ldots,l_r\}$ of a system of parameters of an unmixed module $M$
which guarantee that 
$M/({l_1,\ldots,l_r}) M$ is an unmixed module, too.

Recall that a module is called {\it unmixed} if all its associated prime
ideals have the same height. This can be characterized cohomologically as
follows where we use the convention that a module of negative dimension is
the zero module. The result extends \cite{N-gorliaison}, Lemma 2.11.

\begin{lemma} \label{char-of-unmixedness} Let $M$ be a finitely generated
  $R$-module of dimension $d$.
Then the following conditions are equivalent:
\begin{itemize}
\item[(a)] $M$ is unmixed.
\item[(b)] $\dim \ER^i(M,R) \leq \dim R - 1 - i$ \; if $\dim R - d < i$.
\item[(c)] $\dim R/\fa_i(M) < i$ \; if $i < d$ where $\fa_i(M) = \Ann
  \HH^i(M)$.
\end{itemize}
\end{lemma}

\begin{proof} Choose a complete intersection $\fc = (f_1,\ldots,f_{\dim R -
      d}) \subset \Ann M$. Then we can continue as in the proof of
    \cite{N-hilbert_under_liaison}, Lemma 4.
\end{proof}

We recall the following notion.

\begin{definition} \label{def-M-basis}
  Let $M$ be a finitely generated $R$-module of dimension
  $d$. Let $I \subset R$ be an  ideal such that $\dim M/I M =
  0$. The set $B = \{a_1,\ldots,a_n\} \subset R$ is said to be an {\it
  $M$-basis} of $I$ if its elements form a minimal basis of $I$ and every
  subset of $B$ consisting of $d$ elements is a system of parameters of $M$.
\end{definition}

Fortunately such a basis often exists.

\begin{lemma} \label{M-basis_exists} Let $I \subset R$ be an ideal and let
  $M_1,\ldots,M_s$ be finitely generated $R$-modules such that $\dim M_i/I
  M_i = 0$ for all $i = 1,\ldots,s$. Then it holds:
\begin{itemize}
\item[(a)] (local case) The ideal $I$ has a minimal basis which is an
  $M_i$-basis for all $i = 1,\ldots,s$.
\item[(b)] (graded case) If $I$ is a homogeneous ideal admitting a
  generating set consisting of homogeneous elements of degree $e$ then $I$
  has a
  minimal basis which is formed by homogeneous elements of degree $e$ and
  is an $M_i$-basis for all $i = 1,\ldots,s$.
\end{itemize}
\end{lemma}

\begin{proof} Claim (a) is Proposition 1.9  and (b) is Proposition 3.3 in
  \cite{SV-buch}. Note that the latter can fail if $K$ is a finite field.
\end{proof}

Often an ideal $I$  is called a parameter ideal
for $M$ if  $\dim M/I M = 0$. We want to extend this notion slightly. We
will say that $I$ is a
{\it parameter ideal} for $M$ if $\dim M/I M = \max \{0, \dim M - \mu(I)\}$
where $\mu(I)$ denotes the number of minimal generators of $I$. This means
that a parameter ideal is either generated by a subsystem of a system of
parameters or 
it contains a system of parameters for $M$.

We introduce one more piece of notation. If $M$ is an $R$-module then we
call $M^{sm} = M/\HH^0(M)$ the {\it slight modification} of $M$. If $M =
R/I$ for an ideal $I$ of $R$ then $M^{sm} = R/J$ where
$J \subset R$ is the saturation of the ideal $I$.  More
generally,  if $\ffi: F \to M$ is an epimorphism where $F$ is a free
$R$-module  and we put $E = \ker \ffi$ then
$$
E^{sat} = \bigcup_k (E :_F \fm^k)
$$
is the saturation of $E$ and $M^{sm} \cong F/E^{sat}$. Note that $M$ and
$M^{sm}$ have the same associated sheaf.

The goal of this section is the next result.

\begin{proposition} \label{unmixedness_under_hypersurface_sections}
Let $\{l_1,\ldots,l_s\}$ be a subsystem of a system of parameters for the
finitely generated $R$-module  $M$ such that $I = (l_1,\ldots,l_s) R$ is a
parameter ideal  for all $\ER^i(M,R)$ where $i \neq \dim R - \dim
M$. If $M$ is a graded module we will also assume that all the elements
$l_i$ are homogeneous of degree $e$.
Suppose that
$M^{sm }$ is unmixed. Then $(M/I M)^{sm}$ is unmixed too. Furthermore,
there is a system of generators
$\{h_1,\ldots,h_s\}$ of $I$, which consists  of homogeneous elements of
degree $e$
in the graded case, such that
$(M/(h_1,\ldots,h_i)M)^{sm }$ is unmixed
  for all $i= 0,\ldots,s$.
\end{proposition}

\begin{proof}  We will show the existence of the elements $h_i$
  successively. Suppose we have found minimal generators
  $h_1,\ldots,h_t$ of $I$ where $0 \leq t <s$
  such that   $(M/(h_1,\ldots,h_i)M)^{sm }$ is unmixed for all $i=
  0,\ldots,t$. Choose elements $f_{t+1},\ldots,f_s \in I$ (of degree $e$ in
  the graded case) such that $I = (h_1,\ldots,h_t)R +
  (f_{t+1},\ldots,f_s)R$. Put $N = M/(h_1,\ldots,h_t) M$ and $J =
  (f_{t+1},\ldots,f_s)R$.

Since $(M/(h_1,\ldots,h_i)M)^{sm }$ is unmixed for all $i=
  0,\ldots,t$, $\{h_1,\ldots,h_t\}$ is not just a subsystem of parameters
  for $M$ but even an $M$-filter regular sequence. Hence \cite{NS1},
  Theorem 3.3 yields for all integers $i$
$$
 \fa_i(M) \cdot \ldots \cdot \fa_{i+t}(M) \subset \fa_i(N).
$$
Since $I$ is a parameter ideal for all $\ER^j(M,R)$ with $j \neq \dim R -
\dim M$ we conclude by local duality that we have for all $i < d - t$ where
$d = \dim M$:
\begin{eqnarray*}
\dim R/(\fa_i(N) + J) & = & \dim R/(\fa_i(N) + I) \\
& \leq & \max \{\dim R/(\fa_j(M) +I) \s i \leq j \leq i + t\} \\
& = & \max \{0, \dim R/\fa_j(M) - s  \s i \leq j \leq i + t\} \\
&  \leq & \max \{0, i+t-1-s \}.
\end{eqnarray*}
Here the latter estimate is due to Lemma~\ref{char-of-unmixedness}.

Now we consider the set
$$
P = \{ i \in \ZZ \s i < \dim N = d-t, \; \dim R/\fa_i(N) = i-1 > 0 \}.
$$
We have $\dim N/J N < \dim N$ and, for all $i \in P$, $\dim R/(\fa_i(N) + J)
<\dim R/\fa_i(N)$ because $t < s$. Thus (using arguments as in
Lemma~\ref{M-basis_exists})
there is an element $h \in J$
(homogeneous of degree $e$ in the graded case) which is a parameter for $N$
and all $\ER^{\dim R -i}(N,R)$ with $i \in P$.

Since $N^{sm }$ is unmixed the module $0 :_N h$ has finite length. Thus the
exact sequence
$$
0 \to 0 :_N h \to N \to N/0 :_N h \to 0
$$
implies
$$
\ER^j(N/0 :_N h, R) \cong \ER^j(N, R) \quad \mbox{for all} \; j \neq \dim
R.
$$
Multiplication by $h$ provides the exact sequence
$$
0 \to N/0 :_N h (-e) \to N \to N/h N \to 0.
$$
Using the isomorphisms above we see that it induces exact sequences
\begin{eqnarray*}
\lefteqn{ \ER^{\dim R - i - 1}(N, R)(-e) \stackrel{h}{\longrightarrow}
\ER^{\dim R -
  i - 1}(N, R) \to \ER^{\dim R - i}(N/h N, R) \to \hspace*{5cm} } \\
& & \hspace*{7cm} \ER^{\dim R - i}(N, R)
\stackrel{h}{\longrightarrow}  \ER^{\dim R - i}(N, R).
\end{eqnarray*}
Let $0 < i < \dim N -1 = d-t-1$. Since $N$ is unmixed we have $\dim \ER^{\dim R
  - i}(N, R) \leq i-1$ and $\dim \ER^{\dim R - i - 1}(N, R) \leq i$ where
for the latter estimate equality holds if and only if $i + 1 \in P$. But
then $h$ is a parameter for $\ER^{\dim R - i - 1}(N, R)$ due to our choice
of $h$.  Thus in either
case we have $\dim \ER^{\dim R - i - 1}(N, R)/ h \ER^{\dim R - i - 1}(N, R)
\leq i-1$. Therefore,
using the exact sequence above we obtain
$$
\dim \ER^{\dim R - i}(N/h N, R) \leq i-1 \quad \mbox{for all} \; i \neq
\dim N -1.
$$
Applying the isomorphisms 
$$
\Ext^j_R(N/ h N, R) \cong \Ext^{j-1}_{R/ h R}(N/ h N, R/ h R)(e) 
$$ 
we get that 
$$
\dim \Ext^{\dim R/ h R - i}_{R/ h R}(N/h N, R/h R) \leq \max \{0, i-1\} \quad \mbox{for all} \; i \neq
\dim N -1.
$$
According to Lemma~\ref{char-of-unmixedness} we conclude that $(N/h
N)^{sm }$ is an unmixed module. Thus we may put $h_{t+1} = h$ and we are
done.
\end{proof}


\section{K\"unneth formulas} \label{Kunneth}

Let $R$ be a graded $K$-algebra where $K$ is an arbitrary field. Let $M$
and $N$ be $R$-modules, not necessarily finitely generated. We want to
relate the local cohomology modules of $M \tk N$ to those of $M$ and
$N$ such that the corresponding homomorphisms are $(R \tk R)$-linear. Hence we
cannot use results of Grothendieck because they would give us only $K$-linear
homomorphisms. Instead we  adapt the method of St\"uckrad and
Vogel (cf.\
\cite{SV-buch}) which they used in order to show a K\"unneth formula for the
Segre product of $M$ and $N$.

We begin with some preliminary results where $R_1$ and $R_2$ denote graded
$K$-algebras.

\begin{lemma} \label{intro_of_tau-null} Let $M_i$ and $N_i$ be graded $R_i$
  modules for $i = 1, 2$. Put $R = R_1 \tk R_2$. Then there is a natural
  graded    $R$-homomorphism
$$
\tau_0: \Hom_{R_1}(N_1, M_1) \tk \Hom_{R_2}(N_2, M_2) \to \Hom_R(N_1 \tk
N_2, M_1 \tk M_2).
$$
\end{lemma}

\begin{proof} Let $p_1, p_2$  be integers and let $f_i \in [\Hom_{R_i}(N_i,
  M_i)]_{p_i}$ for $i = 1, 2$. We define for homogeneous elements $n_i \in
  N_i$ where $i = 1, 2$:
$$
\tau_0(f_1 \tk f_2) (n_1 \tk n_2) = f_1(n_1) \tk f_2(n_2).
$$
It is clear that this induces a well-defined $R$-homomorphism
$$
\tau_0(f_1 \tk f_2): N_1 \tk N_2 \to  M_1 \tk M_2
$$
of degree $p_1 + p_2$. Hence it provides a well-defined $K$-linear map
$$
\tau_0: \Hom_{R_1}(N_1, M_1) \tk \Hom_{R_2}(N_2, M_2) \to \Hom_R(N_1 \tk
N_2, M_1 \tk M_2)
$$
which preserves degrees.
Straightforward computations show that $\tau_0$ is even a natural
$R$-homomorphism.
\end{proof}

The map of the lemma is often very nice.

\begin{lemma} \label{tau-null_is_iso} With the notation and assumptions of
  Lemma~\ref{intro_of_tau-null} suppose additionally that $N_i$ is a
  finitely generated $R_i$-module for $i = 1, 2$. Then the map $\tau_0$ is
  an isomorphism.
\end{lemma}

\begin{proof} First, we assume that the modules $N_i$ are free. Thus, since the
  map $\tau_0$ is natural we can even assume that $N_i = R_i(p_i)$ with
  $p_i \in \ZZ$ for $i = 1, 2$. Then the claim is easy to check.

Second, we consider the general case. Let
$$
G_i \to F_i \to N_i \to 0
$$
be free graded presentations. Put $H_i = \Hom_{R_i}(N_i, M_i)$ for $i = 1,
2$. We we have exact sequences
$$
0 \to H_i \to \Hom_{R_i}(F_i, M_i) \to \Hom_{R_i}(G_i, M_i) \quad (i = 1,
2)
$$
and
$$
(G_1 \tk F_2) \oplus (G_2 \tk F_1) \to F_1 \tk F_2  \to N_1 \tk N_2 \to 0.
$$
Since the functor \underline{ }\,\,$\tk$\,\,\underline{ } is exact they
induce the following commutative diagram with exact rows
$$
\begin{array}{cccccc}
0 \to &  H_1 \tk H_2 & \to & \Hom_{R_1}(F_1, M_1) \tk \Hom_{R_2}(F_2, M_2) &
\hspace*{1cm}\\
& \downarrow \tau_0 & & \downarrow \bar{\tau}_0 & \\
0 \to & \Hom_R(N_1 \tk N_2, M_1 \tk M_2) & \to & \Hom_R(F_1 \tk F_2, M_1
\tk M_2) &
\end{array}
$$
$$
\begin{array}{cccc}
\hspace*{1cm} &\to &  (\Hom_{R_1}(F_1, M_1) \tk \Hom_{R_2}(G_2, M_2)) \oplus
(\Hom_{R_2}(F_2, M_2) \tk \Hom_{R_1}(G_1, M_1)) \\
& & \downarrow \hat{\tau}_0 \oplus \tilde{\tau}_0 \\
& \to & \Hom_R((G_1 \tk F_2) \oplus (G_2 \tk F_1), M_1 \tk M_2)
\end{array}
$$
where $\bar{\tau}_0, \hat{\tau}_0, \tilde{\tau}_0$ are the corresponding
natural homomorphisms. These maps are even isomorphisms by the first part
of the proof. Thus $\tau_0$ is an isomorphism as well.
\end{proof}

For a graded  $R$-module $M$ we denote its graded injective hull by $E(M)$.
As the last preparation we need.

\begin{lemma} \label{product_of_injectiv_modules} Let $I_i$ be graded
  injective $R_i$-modules for $i = 1, 2$. Let $\fm$ be the irrelevant
  maximal ideal of $R = R_1 \tk R_2$. Then we have
$$
\HH^j(I_1 \tk I_2) = 0 \quad \mbox{for all} \; j \geq 1.
$$
\end{lemma}

\begin{proof}  According to results of Matlis \cite{Matlis_injective} we
  know that $I_i$ is a direct sum of injective hulls $E(R_i/\fp_i)(p_i)$
  where $\fp_i \in \Proj (R_i) \cup \{\fm_i\}$ and $p_i \in \ZZ$. Since the
  tensor product and local cohomology commute with direct sums we may
  assume that $I_i = E(R_i/\fp_i)(p_i)$.

If $\fp_i = \fm_i$ for $i = 1, 2$
  then we have $\Supp I_1 \tk I_2 = \{\fm\}$. This implies $\HH^j(I_1 \tk
  I_2) = 0$ for all $j \geq 1$.

Otherwise we can find homogeneous elements $y_i \in [R_i]_{p_i} \backslash
\fp_i$ such that
$y = y_1 \tk y_2$ has degree $p > 0$, i.e.\ $y \in \fm$. Due to our choice
of the elements $y_i$ the multiplication gives isomorphisms $I_i
\stackrel{y_i}{\longrightarrow} I_i(p_i)$. Therefore the multiplication $I
\stackrel{y}{\longrightarrow} I(p)$ where $I = I_1 \tk I_2$ is an isomorphism
as well. Thus we obtain for all $j \geq 1$ isomorphisms $\HH^j(I)
\stackrel{y}{\longrightarrow} \HH^j(I)(p)$. Since $\Supp \HH^j(I) \subset
\{\fm\}$ it follows that $\HH^j(I) = 0$ for all $j \geq 0$.
\end{proof}

Now we are in the position to state and to prove the main result of this
section.

\begin{theorem} \label{Kuenneth_formula-thm} Let $M_i$ and $N_i$ be graded
  $R_i$-modules and let $R = R_1 \tk R_2$. Then we have for all $k \geq 0$:
  \begin{itemize}
\item[(a)] There are natural graded $R$-homomorphisms
$$
\tau_k: \bigoplus_{i+j = k} \Ext^i_{R_1}(N_1, M_1) \tk \Ext^j_{R_2}(N_2, M_2)
\to \Ext^k_R(N_1 \tk N_2, M_1 \tk M_2).
$$
\item[(b)] (K\"unneth formulas)  There are natural graded $R$-isomorphisms
$$
\sigma_k: \bigoplus_{i+j = k} H^i_{\fm_1}(M_1) \tk H^j_{\fm_2}(M_2) \to
\HH^k(M_1 \tk M_2).
$$
\end{itemize}
\end{theorem}

\begin{proof} First we prove (a). Let $I^{\bullet}_i$, $i= 1, 2$, be
  injective resolutions of
  $M_i$. Put $P^{\bullet}_i = \Hom_{R_i}(N_i, I^{\bullet}_i)$.  For the
  (co)homology modules of these complexes we have by
  \cite{Cartan-Eilenberg}, Theorem~IV.7.2
$$
\bigoplus_{i+j = k} \Ext^i_{R_1}(N_1, M_1) \tk \Ext^j_{R_2}(N_2, M_2) =
\bigoplus_{i+j = k} H^i(P^{\bullet}_1) \tk H^J(P^{\bullet}_2 \cong
H^k(P^{\bullet}_1 \tk P^{\bullet}_2).
$$
According to Lemma~\ref{intro_of_tau-null} we have natural homomorphisms of
complexes of $R$-modules
$$
P^{\bullet}_1 \tk P^{\bullet}_2 = \Hom_{R_1}(N_1, I^{\bullet}_1) \tk
\Hom_{R_2}(N_2, I^{\bullet}_2) \to \Hom_R(N_1 \tk N_2, I^{\bullet}_1 \tk
I^{\bullet}_2).
$$
Since $I^{\bullet}_1 \tk I^{\bullet}_2$ is an acyclic right complex with
$H^0 (I^{\bullet}_1 \tk I^{\bullet}_2) = M_1 \tk M_2$ there are
homomorphisms
$$
H^k (P^{\bullet}_1 \tk P^{\bullet}_2) \to H^k (\Hom_R(N_1 \tk N_2,
I^{\bullet}_1 \tk I^{\bullet}_2) \to \Ext_R^k(N_1 \tk N_2, M_1 \tk M_2)
$$
(cf., for example, \cite{Cartan-Eilenberg}, Proposition V.1.1a). Now
composition of the maps described above gives the desired homomorphisms.

Second we prove (b). Let $s$ be an integer such that the ideals $\fm_1$ and
$\fm_2$ can be generated by less than $s$ elements. Then we have for all
positive integers $t$ canonical epimorphisms
$$
R/\fm^{s t} \to R_1/\fm_1^t \tk R_2/\fm_2^t \to R/\fm^t.
$$
They induce the homomorphisms
$$
\Ext^k_R(R/\fm^{s t}, M_1 \tk M_2)  \to \Ext^k_R(R_1/\fm_1^t \tk
R_2/\fm_2^t, M_1 \tk M_2) \to \Ext^k_R(R/\fm^t, M_1 \tk M_2).
$$
Taking the direct limit over all $t \geq 1$ we obtain
$$
\lim_{\longrightarrow} \Ext^k_R(R_1/\fm_1^t \tk R_2/\fm_2^t, M_1 \tk M_2) \cong
\HH^k(M_1 \tk
M_2).
$$
Now we use part (a). It provides natural homomorphisms
$$
\Ext^k_{R_1}(R_1/\fm_1^t, M_1) \tk \Ext^k_{R_2}(R_2/\fm_2^t, M_2) \to
\Ext^k_R(R_1/\fm_1^t \tk R_2/\fm_2^t, M_1 \tk M_2).
$$
Taking again the direct limit and using the isomorphism above we get
natural homomorphisms
$$
\sigma_k: \bigoplus_{i+j = k} H^i_{\fm_1}(M_1) \tk H^j_{\fm_2}(M_2) \to
\HH^k(M_1 \tk M_2).
$$
According to \cite{Cartan-Eilenberg}, Proposition V.4.4,  claim (b) holds,  if
$\sigma_0$ is an isomorphism for all modules $M_1$ and $M_2$ and $\sigma_k$
is an isomorphism for all injective modules $M_1$ and $M_2$. The second
statement is true by Lemma~\ref{product_of_injectiv_modules}. Thus it
remains to show that $\sigma_0$ is an isomorphism.

Due to Lemma~\ref{tau-null_is_iso} the maps
$$
\tau_0: \Hom_{R_1}(R_1/\fm_1^t, M_1) \tk \Hom_{R_2}(R_2/\fm_2^t, M_2) \to
\Hom_R(R_1/\fm_1^t \tk R_2/\fm_2^t, M_1 \tk M_2)
$$
are isomorphisms. Thus taking the direct limit shows that $\sigma_0$ is an
isomorphism.
\end{proof}

The last result allows us to relate information about two projective
subschemes to properties of its embedded join. We state the results in
greater generality. First we need.

\begin{lemma} \label{CMness_of_join} Let $M_i$ be graded $R_i$-modules for
  $i = 1, 2$. Then we have
\begin{itemize}
\item[(a)] $M_1 \tk M_2$ is a Cohen-Macaulay module if and only if $M_1$
  and $M_2$ are Cohen-Macaulay.
\item[(b)] $R = R_1 \tk R_2$ is a Gorenstein algebra if and only if $R_1$
  and $R_2$ are Gorenstein.
\end{itemize}
\end{lemma}

\begin{proof} Write $R_i = S_i/I_i$ where $S_i$ is a polynomial ring over
  $K$ and
  $I_i \subset S_i$ is a homogeneous ideal. Let $F_{\bullet}^{(i)}$
  denote the minimal free resolution of $M_i$ as $S_i$-module. Then
  $F_{\bullet}^{(1)} \tk F_{\bullet}^{(2)}$ is a minimal free resolution of
  $M_1 \tk M_2$ as $S_1 \tk S_2$-module. In conjunction with the
  Auslander-Buchsbaum
  formula claim (a) follows. Moreover, we see that the Cohen-Macaulay type
  of $M_1 \tk M_2$ is the product of the Cohen-Macaulay types of $M_1$ and
  $M_2$. Hence claim (b) follows from the fact that $R$ is Gorenstein if
  and only if its Cohen-Macaulay type is one.
\end{proof}

Now we can show.

\begin{proposition} \label{join-properties}  Let $M_i$ be graded
  $R_i$-modules for
  $i = 1, 2$. Assume that the rings $R_i$ are Gorenstein. Then we have
\begin{itemize}
\item[(a)] $M_1 \tk M_2$ is an unmixed $R_1 \tk R_2$-module if and only if
  $M_1$ and $M_2$ are unmixed.
\item[(b)] $\depth M_1 \tk M_2 = \depth M_1 + \depth M_2$.
\item[(c)] The following conditions are equivalent:
\begin{itemize}
\item[(i)] $M_1 \tk M_2$ is unmixed and locally Cohen-Macaulay.
\item[(ii)]  $M_1 \tk M_2$ is  Cohen-Macaulay.
\item[(iii)]  $M_1$ and $M_2$ are Cohen-Macaulay.
\end{itemize}
\end{itemize}
\end{proposition}

\begin{proof} (a) Since  $R = R_1 \tk R_2$ is Gorenstein by the previous lemma
  we can apply Lemma~\ref{char-of-unmixedness}. Now suppose that $M_1$ is
  not unmixed. Then there is an $i > \dim R_1 -\dim M_1$ such that $\dim
  \Ext^i_{R_1}(M_1, R_1) \geq \dim R_1 - i$. Let $j = \dim R_2 -\dim
  M_2$. Then we have $\dim \Ext^j_{R_2}(M_2, R_2) = \dim M_2$ and thus
$$
\dim \Ext^i_{R_1}(M_1, R_1) \tk \Ext^j_{R_2}(M_2, R_2) \geq \dim R_1 - i +
\dim M_2 = \dim R - i -j.
$$
Hence Theorem~\ref{Kuenneth_formula-thm} implies
$$
\dim \Ext^{i+j}_R(R, M_1 \tk M_2) \geq \dim R - i -j.
$$
Since $i+j > \dim R - \dim M_1 \tk M_2$, Lemma~\ref{char-of-unmixedness}
shows that $M_1 \tk M_2$ is not unmixed.

Conversely, if $M_1$ and $M_2$ are unmixed then $M_1 \tk M_2$ is unmixed
too using similar arguments as above.

(b) follows from the cohomological characterization of depth and the
K\"unneth formula.

(c) According to Lemma~\ref{CMness_of_join} it suffices to show that (i)
implies (iii). Suppose that $M_1$ is not Cohen-Macaulay. Then there is an
$i < \dim M_1$ such that $\HH^i(M_1) \neq 0$. Let $j = \dim M_2$. Since
$H^j_{\fm_1}(M_1)$ is not finitely generated, the module $\HH^{i+j}(M_1 \tk
M_2)$ is not a finitely generated $R$-module and,  in particular, is not of
finite length. Hence $M_1 \tk M_2$ cannot be equidimensional and locally
Cohen-Macaulay which completes the proof.
\end{proof}

In the special case where $M_1$ and $M_2$ are $K$-algebras the conclusion
from properties of $M_1$ and $M_2$ to properties of $M_1 \tk M_2$ can also
be found in \cite{Vogel-Tata}, Proposition 1.47. 

We remark that the results of this section hold also true for modules over
local rings containing a field $K$. 


\section{A Bezout theorem for modules} \label{Bezout}

In this section $R$ will always denote the polynomial ring $K[x_0,\ldots,x_n]$
with its standard grading or a regular local ring containing a field.
Moreover, $M$ and $N$ will denote finitely 
generated graded $R$ modules and $\Delta$ the diagonal ideal of $R \tk R$. 

We will apply the results of Section \ref{sec-unmixedness} where we
assumed that the ground field $K$ is infinite. This assumption is harmless
because  the Cohen-Macaulay property and Hilbert functions are  not
effected by field extension. Thus the results of the rest of the paper hold
true for arbitrary fields. 

Recall that there is an inequality
$$
\dim M \tr N \leq \dim M + \dim N -\dim R.
$$
If equality holds then it is said that $M$ and $N$ intersect properly. We
need an even stronger condition.

\begin{definition} \label{very_proper_inter-def} It is said that $M$ and
  $N$ intersect {\it very properly} if
$$
\dim M + \dim N \geq \dim R
$$
and
$$
\dim \ER^i(M, R) \tr \ER^j(N, R) = \max \{0, \dim \ER^i(M, R) + \dim
\ER^j(N, R) - \dim R \}
$$
for all integers $i$ and $j$.
\end{definition} 

Similarly, we say that subschemes $X, Y \subset \PP^n$  intersect very
properly if the corresponding modules $R/I_X$ and $R/I_Y$  meet very properly. 

We want to compare the two conditions proper and very proper
intersection. For this we recall that the canonical 
module of $M$ is
$$
K_M = \ER^{n+1 - \dim M}(M, R)(-n).
$$
The sets of top-dimensional associated prime ideals of $M$ and $K_M$ coincide.

\begin{lemma} \label{proper_versus_very-proper_lemma} It holds:
\begin{itemize}
\item[(a)] If $M$ and $N$ intersect very properly then they intersect
  properly.
\item[(b)] If $M$ and $N$ are unmixed locally Cohen-Macaulay modules which
  intersect properly then they intersect very properly.
\end{itemize}
\end{lemma}

\begin{proof} (a) The assumption implies
$$
\dim M + \dim N - (n+1) = \dim K_M \tr K_N.
$$
Hence $\dim K_M \tr K_N = \dim M \tr N$ proves the claim.  \\
(b) Since $M$ and $M$ meet properly we have
$$
\dim M + \dim N - (n+1) = \dim M \tr N = \dim K_M \tr K_N \geq 0.
$$
If $(i, j) \neq (n+1 - \dim M, n+1 - \dim N)$ then one of the factors of
$$
\ER^i(M, R) \tr \ER^j(N, R)
$$
has finite length because of the assumptions on $M$ and $N$. It follows
that
$$
\dim \left ( \ER^i(M, R) \tr \ER^j(N, R) \right ) = 0.
$$
This shows that $M$ and $N$ intersect very properly.
\end{proof}

We will need a result about the change of the degree under hyperplane
section. Let us denote the dimension of the $R$-module $M$ by $e > 0$. Then we
can write its Hilbert polynomial as
$$
p_M(j) = h_0 \cdot \binom{j+e}{j} + h_{1} \cdot \binom{j+e-1}{j} + \ldots +
h_e \cdot \binom{t}{t}
$$
with integers $h_0,\ldots,h_{e}$.  The degree of $M$ is
$$
\deg M = \left \{ \begin{array}{ll}
h_0 & \mif \dim M > 0 \\
length (M) & \mbox{otherwise}.
\end{array} \right.
$$

Since we are not aware of a reference for the following observation in the
generality we need, we state it explicitly. 

\begin{lemma} \label{degree_and_hypersurface-section} Let $f \in R$ be a
  homogeneous parameter for $M$ of degree $k$. Then we have
$$
\deg M/f M  = \left \{ \begin{array}{ll}
k \cdot \deg M & \mif f \notin \fp \; \mbox{for all} \; \fp \in \Ass M \\
& \mbox{where} \;  \dim R/\fp = \dim M -1 \\
k \cdot \deg M + \deg (0 :_M f) & \mbox{otherwise}.
\end{array} \right.
$$
\end{lemma}

\begin{proof} The exact sequence
$$
0  \to (0 :_M f)(-k) \to M(-k) \stackrel{f}{\longrightarrow} M \to M/f M
\to 0
$$
provides for the Hilbert polynomials
$$
p_{M/f M}(j) = p_M(j) - p_M(j-k) + p_{0 :_M f}(j-k).
$$
Since $\dim M/f M > \dim (0 :_M f)$ if and only if $f \notin \fp$ for all $\fp
\in \Ass M$ where $\dim R/\fp = \dim M -1$ the claim follows by comparing
the coefficients of the polynomials above.
\end{proof}

We will refer to the next result as a Bezout theorem. The name will become
clear from the consequences of the result for the intersection of
projective schemes. In the proof we will use the isomorphism
$$
M \tr N \cong (M \tk N)/ \Delta (M \tk N)
$$
frequently. Its use is often called diagonal trick or reduction to the
diagonal. Analogous to the case of subschemes we call $M \tk N$ the join of
$M$ and $N$ and $M \tr N$ their intersection.

\begin{theorem} \label{Bezout_for_modules} Suppose that $M$ and $N$ are
  unmixed modules which intersect very properly. Then $(M \otimes_R
  N)^{sm}$ is an unmixed module.

Moreover, if $M \otimes_R N$ has positive
  dimension or $\depth M + \depth N \geq n+1$ then we have
$$
\deg M \otimes_R N = \deg M \cdot \deg N.
$$
\end{theorem}

\begin{proof} Let $S = R \tk R$. Then we have by
  Theorem~\ref{Kuenneth_formula-thm} the isomorphisms
\begin{eqnarray*}
\lefteqn{\bigoplus_{i+j = k} \ER^i(M, R) \tr \ER^j(N, R)} \\
 & \cong & (\bigoplus_{i+j =
  k} \ER^i(M, R) \tk \ER^j(N, R))/ \Delta (\bigoplus_{i+j =
  k} \ER^i(M, R) \tk \ER^j(N, R)) \\
& \cong & \Ext^k_S(M \tk N, S)/ \Delta \Ext^k_S(M \tk N, S).
\end{eqnarray*}
Since $\Delta$ has $n+1$ minimal generators and $M$  and $N$ intersect very
properly we see that $\Delta$ is a
parameter ideal for all $\Ext^k_S(M \tk N, S)$ and for $M \tk N$. Therefore
Proposition~\ref{unmixedness_under_hypersurface_sections} shows that
$$
((M \tk N)/ \Delta (M \tk N))^{sm} \cong (M \tr N)^{sm}
$$
is an unmixed module proving the first claim. \\
In order to show the claim on the degrees we show firstly
$$
\deg M \tk N = \deg M \cdot \deg N.
$$
Let $e$ and $f$ denote the dimension of $M$ and $N$, respectively.
For all integers $t \geq 0$ we have the following relation of Hilbert
functions:
$$
h_{ M \tk N}(t) = \sum_{i+j=t} h_M(i) \cdot h_N(j).
$$
Let $r$ be an integer such that the Hilbert functions of $M, N$ in
degree $i > r$ are given by their respective Hilbert polynomials. Then we
obtain for $t \geq 2r$:
\begin{eqnarray*}
h_{ M \tk N}(t)  & = &  \sum_{i = 0}^t p_M(i) \cdot p_N(t-i) + \sum_{i=0}^r
[h_M(i)-p_M(i)] \cdot p_N(t-i) \\
& & \mbox{} + \sum_{i= t-r}^t p_M(i) \cdot [h_N(t-i) -
p_N(t-i)] \\
& = & \deg M \cdot \deg N \cdot \sum_{i=0}^t \binom{i}{e-1} \cdot
\binom{t-i}{f-1} + O(t^{e+f-2}) \\
& = & \deg M \cdot \deg N \cdot \binom{t}{e+f-1} + O(t^{e+f-2}).
\end{eqnarray*}
It follows that $\deg M \tk N = \deg M \cdot \deg N$ as claimed.

Above we have already  applied
Proposition~\ref{unmixedness_under_hypersurface_sections}. Now we use its
full strength. It provides a minimal basis $\{h_0,\dots,h_n\}$ of $\Delta$
such that
$$
((M \tk N)/(h_0,\ldots,h_i) (M \tk N))^{sm}
$$
is an unmixed module for all $i = 0,\ldots,n$. This implies
$$
\deg M \tk N = \deg (M \tk N)/(h_0,\ldots,h_n)(M \tk N) = \deg M \tr N
$$
according to our assumptions and the previous lemma.
\end{proof}

Note that the degree relation in the last statement is very much in the
spirit of Bezout's original result on the intersection of plane
curves. Nowadays it is well-known that a similar formula cannot be true for
arbitrary proper intersections. In general one has to replace the length
multiplicity by a suitable intersection multiplicity   (cf.\
\cite{Serre-algebre-locale}, \cite{Vogel-Tata}
\cite{Kirby_Bezout}, Theorem 2.8, \cite{Fulton-inter_book} and
\cite{Roberts-book}). 
However, our result gives a condition when the simplest multiplicity works
and in that case the extra information that the intersection is unmixed up
to an irrelevant component.   

\begin{example} \label{example-skew_lines} Let $X \subset \PP^3$ be the
  union of two skew lines defined by the ideal $I = (x_0, x_1) \cap (x_2,
  x_3)$ of the polynomial ring $R = K[x_0,\ldots,x_3]$ and let $Y \subset
  \PP^3$ be the complete intersection defined by $J = (x_1 + x_2, x_0 +
  x_3)$. Then $M = R/I$ and $N = R/J$ intersect very properly in dimension
  zero. However, it holds
$$
\deg M \tr N = 3 \neq 2 = \deg M \cdot \deg N.
$$
Hence this  example  shows that in general the degree relation does not
hold true in the
preceding statement if $\dim M \tr N = 0$.
\end{example}

The previous result implies for subschemes of $\PP^n$.

\begin{corollary} \label{Bezout_in_Pn-corollary} Let $X$ and $Y$ denote
  equidimensional subschemes of $\PP^n$ which intersect very properly. Then
  $X \cap Y$ is an equidimensional subscheme.

Moreover, if $X \cap Y \neq \emptyset$ then we have
$$
\deg X \cap Y = \deg X \cdot \deg Y.
$$
\end{corollary}


A special case of Bertini's theorem says that a sufficiently general
hyperplane section of an equidimensional subscheme is again
equidimensional. 
The first claim generalizes this statement and makes  the
condition  ``sufficiently general'' more precise.   
It  implies in particular that $X \cap Y$ does not have embedded components.

Taking Proposition~\ref{join-properties} into account the proof of
Theorem~\ref{Bezout_for_modules} shows.

\begin{proposition} \label{depth_of_intersection_lemma-prop} Let $M$ and $N$ be
  unmixed modules which intersect very properly. Then it holds
$$
\depth M \tr N = \max \{0, \depth M + \depth N - (n+1) \}.
$$
\end{proposition}

If $\depth M + \depth N \geq n+1$  the proposition implies
$$
\depth M + \depth N = \depth R + \depth M \tr N,
$$
i.e., the modules $M$ and $N$ satisfy the depth formula in the sense of
Huneke and Wiegand \cite{Huneke-Wiegand-Math.Scand}.


\section{Lifting the Cohen-Macaulay property} \label{lifting}

Here we look for conditions/results which say essentially that the
Cohen-Macaulay 
property of the very proper intersection of $M$ and $N$ forces $M$ and $N$
to be Cohen-Macaulay modules. We use the same notation as in the previous
section. 

A first result follows immediately from the diagonal trick. In fact, if $M$
and $N$ intersect properly then any minimal generating set of the diagonal
ideal 
is a subsystem of a system of parameters for $M \tk N$. Thus the
isomorphism $M \tr N \cong  (M \tk N)/ \Delta (M \tk N)$ implies that $M
\tk N$ is Cohen-Macaulay. Thus we have seen. 

\begin{lemma} \label{lemma-Serre-lifting} Let $M$ and $N$ denote 
  modules which intersect  properly. If $M \tr N$ is Cohen-Macaulay of
  positive dimension then $M$ and $N$ are Cohen-Macaulay. 
\end{lemma} 

\begin{remark} The result implies for subschemes $X, Y \subset \PP^n$: If they meet
properly in a non-empty scheme and $R/(I_X + I_Y)$ is Cohen-Macaulay then
$X$ and $Y$ are arithmetically Cohen-Macaulay. For a result where we only
assume that $X \cap Y$
is arithmetically Cohen-Macaulay, we refer to Corollary
\ref{CM-lifting-schemes-cor}. 

Note that in case $R/(I_X + I_Y)$ is Cohen-Macaulay,  the scheme $X \cap Y$
is arithmetically Cohen-Macaulay. In general, the converse is not true
because the ideal $I_X + I_Y$ is not necessarily saturated. 
\end{remark} 

In order to generalize the statement above we begin with the following
observation. 

\begin{lemma} \label{Huneke-U-lifting} Let $M$ denote an unmixed graded
  $R$-module
and let $f \in R$ be a homogeneous $M$-regular element. Suppose that $(M/f
M)^{sm}$ has depth $t \geq 2$. Then we have $\depth M = t+1$ and $M/f M = (M/f
M)^{sm}$.

In particular, if $(M/fM)^{sm}$ is a Cohen-Macaulay module of dimension
$\geq 2$ then $M$ is Cohen-Macaulay.
\end{lemma}

\begin{proof} In the special  case where $M$ is a ring our claim  follows
  from \cite{HU-hyperplane-sect}, Proposition 2.1. The proof in the general
  case is similar to those of Huneke and Ulrich. Since it is short we give
  it for the reader's convenience. 

By the definition of a slight modification we have an exact sequence
$$
0 \to \HH^0(M/fM) \to M/fM \to (M/fM)^{sm} \to 0
$$
which implies $\HH^i(M/fM) \cong \HH^i((M/fM)^{sm})$ for all $i \geq
1$. In particular, it follows $\HH^1(M/f M) = 0$.

Let $k$ denote the degree of $f$. Then the exact sequence
$$
0 \to M(-k) \stackrel{f}{\longrightarrow} M \to M/f M \to 0
$$
provides  the  exact cohomology sequence
$$
0 \to \HH^0(M/f M) \to \HH^1(M)(-k) \stackrel{f}{\longrightarrow} \HH^1(M) \to
\HH^1(M/f M) \to \ldots.
$$
According to Lemma~\ref{char-of-unmixedness} the module $\HH^1(M)$ is
finitely generated. Since $\HH^1(M/f M) = 0$, Nakayama's lemma yields
$\HH^1(M) = 0$. It follows that  $\HH^0(M/f M) = 0$ which implies our claims.
\end{proof}

Now we are ready for the main result of this section. It may be viewed as a
vast generalization of the previous lemma. However, the lemma is used in
the proof below.

\begin{theorem} \label{CM-lifting_theorem} Let $M$ and $N$ denote unmixed
  modules which intersect very properly. Suppose that one of the following
  conditions is satisfied:
\begin{itemize}
\item[(i)] $(M \tr N)^{sm}$ is Cohen-Macaulay and has  dimension $\geq 2$.
\item[(ii)] $\dim M \tr N = 0$  and $\deg M \tr N = \deg M
  \cdot \deg N$.
\end{itemize}
Then $M$ and $N$ are Cohen-Macaulay modules.
\end{theorem}

\begin{proof} Since $M$ and $N$ intersect very properly we can apply
  Proposition~\ref{unmixedness_under_hypersurface_sections} as in the proof
  of Theorem~\ref{Bezout_for_modules}. Thus there is a minimal basis
  $\{h_0,\ldots,h_n\}$ of the diagonal ideal $\Delta$ consisting of linear
  forms such that
$$
((M \tk N)/(h_0,\ldots,h_i) (M \tk N))^{sm}
$$
is an unmixed module for all $i = 0,\ldots,n$. Moreover we have
$$
(M \tk N)/(h_0,\ldots,h_n) (M \tk N) \cong M \tr N.
$$

Now suppose that assumption (ii) is fulfilled. Then $(M \tk
N)/(h_0,\ldots,h_{n-1}) (M \tk N)$ is a module of dimension one and degree
$\deg M \cdot \deg N$. Since the module $M \tr N$ has the same degree
we conclude by Lemma~\ref{degree_and_hypersurface-section} that the
irrelevant maximal ideal is not an
associated prime ideal of  $(M \tk N)/(h_0,\ldots,h_{n-1}) (M \tk N)$,
i.e., $(M \tk N)/(h_0,\ldots,h_{n-1}) (M \tk N)$ is Cohen-Macaulay. It
follows that $M \tk N$ is Cohen-Macaulay. 

Next, assume that condition (i) is satisfied. Then we can apply
Lemma~\ref{Huneke-U-lifting} successively. It follows that $(M \tk N)^{sm}$
is Cohen-Macaulay. Due to Proposition~\ref{CMness_of_join}(a) we know that $M
\tk N$ is an unmixed module, thus $(M \tk N)^{sm} = M \tk N$.

Therefore we have seen that any of the conditions (i) and  (ii) implies that
$M \tk N$ is a Cohen-Macaulay module. Thus
Lemma~\ref{CMness_of_join}(a) proves our assertion.
\end{proof}

The next result was one of the driving forces of this note.

\begin{corollary} \label{CM-lifting-schemes-cor} Let $X, Y \subset \PP^n$
  denote equidimensional subschemes which meet very properly in an
  arithmetically Cohen-Macaulay subscheme whose dimension is at least one. Then
  $X$ and $Y$ are arithmetically Cohen-Macaulay.

Furthermore, if the Cohen-Macaulay type of $X \cap Y$ is a prime number
then $X$ or $Y$ is arithmetically Gorenstein. If $X \cap Y$ is
arithmetically  Gorenstein then $X$ and $Y$ are arithmetically Gorenstein.
\end{corollary}

\begin{proof} With regard to the previous theorem it suffices to note that
  the Cohen-Macaulay type of $X \cap Y$ is the product of the types of $X$
  and $Y$.
\end{proof}

\begin{remark} \label{CM-lifting-fails-remark} The theorem and its
  corollary do  not hold
  true  when we replace very proper intersection by the weaker
  assumption of a proper intersection. Consider the following example: Let
  $\hat{X}, \hat{Y} \subset \PP^5$ be the cones over the corresponding
  subschemes $X, Y \subset \PP^3$ in Example~\ref{example-skew_lines}. Then
  $\hat{X}$ is locally Cohen-Macaulay at all points except at its vertex. Thus
  $\hat{X}$ and $\hat{Y}$ meet properly but not very properly.  The
  intersection is  an arithmetically Cohen-Macaulay curve but $X$ is not
  even locally Cohen-Macaulay.
\end{remark}


\section{Some splitting criteria for reflexive sheaves} \label{splitting}

In this section we apply the methods, which we have developed in the previous
sections, in order to study tensor products of reflexive sheaves on projective
 space.

To start with we note that two maximal $R$-modules $M, N$ intersect always
 properly. According to Lemma \ref{proper_versus_very-proper_lemma} $M$ and
$N$
meet even very properly if $M$ and $N$ are locally free.

Recall that the module $H^0(M)$ is defined as
$$
H^0(M) = \lim_{\longrightarrow} \Hom_R(\fm^t, M).
$$
It is again a graded module which fits into the exact sequence
$$
0 \to \HH^0(M) \to M \to H^0(M) \to \HH^1(M) \to 0.
$$
In particular it holds (cf., for example, \cite{SV-buch}, Lemma 0.1.8)
$$\HH^i(H^0(M)) \cong \left \{\begin{array}{ll}
\HH^i(M) & \mif i \geq 2 \\
0 & \mif i \leq 1.
\end{array} \right.
$$

The next result shows essentially that the tensor product has nice
properties if the two modules intersect very properly. 

\begin{proposition} \label{proper_makes_tensor_product_nice} Let $M$ and $N$
be maximal $R$-modules which intersect very properly. Then we have:
\begin{itemize}
\item[(a)] If $M$ and $N$ are torsion-free then $(M \tr N)^{sm}$ is
torsion-free.
\item[(b)] If $M$ and $N$ are reflexive then $H^0(M \tr N)$ is reflexive.
\end{itemize}
\end{proposition}

 \begin{proof} Claim (a) is a special case of Theorem \ref{Bezout_for_modules}
because a maximal $R$-module is torsion-free if and only if it is unmixed. 

In order to prove claim (b) we recall that a maximal $R$-module $E$ is
reflexive
if and only if $\dim R/\fa_i(E) \leq i-2$ for all $i < \dim R = n+1$ (cf.\
\cite{EG-buch}, Theorem 3.6). Thus, the K\"unneth formulas provide that the
join
$M \tk N$ is a reflexive $(R \tk R)$-module. Since $M$ and $N$ intersect very
properly the diagonal $\Delta$ is a parameter ideal for $M \tk N$ and all
$\Ext^i_{R \tk R}(M \tk N, R \tk R)$ where $i \neq 0$. Thus, arguing as in
the
proof of Proposition \ref{join-properties} we obtain
$$
\dim \Ext^{n+1-i}_{R}(M \tr N, R) \leq \max \{0, i - 2\} \quad
\mbox{for all} \; i \neq 0.
$$
Hence $H^0(M \tr N)$ is a reflexive $R$-module.
\end{proof}

In order to discuss the last result we consider some examples where we used
the
computer algebra program MACAULAY \cite{macaulay-comp} for carrying out the
computations. 

\begin{example} \label{examples_reflexive_sheaves}
Let $I_1$ and $I_2$ denote the homogeneous ideals in
$R = K[x_0,\ldots,x_3]$ of two different points in $\PP^3$. Let $E_1$ and $E_2$ be
the first syzygy modules of $I_1$ and $I_2$, respectively.

Then $I_1 \tr I_1$ is a proper but not a very proper intersection. It turns
out
that $I_1 \tr I_1 \cong (I_1 \tr I_1)^{sm}$ is not torsion-free.

Now let us consider $I_1 \tr I_2$. It is a very proper intersection.
However, it
is not torsion-free. On the other hand $(I_1 \tr I_2)^{sm}$ is torsion-free,
but not reflexive, whereas $H^0(I_1 \tr I_2)$ is even a reflexive $R$-module.

Looking at the reflexive modules $E_1$ and $E_2$ we observe that $E_1 \tr E_1$
is a proper but not a very proper intersection.  Furthermore, $E_1 \tr E_1
\cong
(E_1 \tr E_1)^{sm}$ is not reflexive. The situation is different for
$E_1 \tr E_2$. It is a very proper intersection and $E_1 \tr E_2 \cong
H^0(E_1 \tr E_2)$ is a reflexive $R$-module.
\end{example}

\begin{remark} \label{remark_very_proper_required}
 The example above shows that the conclusions in Proposition
\ref{proper_makes_tensor_product_nice} are not true if we don't assume that
the intersections are very proper.  
\end{remark} 

Let $\cF$ be a sheaf on $\PP^n$. As usual we put for all integers $i$ 
$$
H^i_*(\cF) = \oplus_{t \in \ZZ} H^i(\PP^n, \cF(t)).  
$$
We say that sheaves $\cF$ and 
$\cG$ on $\PP^n$ intersect very properly if the modules  $H^0_*(\cF)$  and
$H^0_*(\cG)$ meet very properly.  Then Proposition
\ref{proper_makes_tensor_product_nice} implies. 

\begin{corollary} \label{cor_tensor_sheaves} Let $\cE,\cF$ denote coherent
  sheaves on $\PP^n$ which intersect very properly. Then we have: 
\begin{itemize} 
\item[(a)] If $\cF$ and $\cG$ are torsion-free then $\cF \otimes \cG$ is
torsion-free. 
\item[(b)] If $\cF$ and $\cG$ are reflexive then $\cF \otimes \cG$ is
  reflexive. 
\end{itemize} 
\end{corollary} 

\begin{remark} \label{remark_refl_sheaves}
 (i) Again, Example \ref{examples_reflexive_sheaves} shows that the
 statements  are no longer true if the intersections are not very proper. \\
(ii) Auslander has shown in \cite{Auslander_Ill_J-1961}, Lemma 3.1 that the
torsion-freeness of $M \tr N$  implies that  $M$ and $N$ are torsion-free.
Huneke and Wiegand (\cite{Huneke-Wiegand-Math.Ann}, Theorem 2.7) proved that
$M$ and $N$ are reflexive if $M \tr N$ is reflexive.

A corresponding statement  does not hold for sheaves. In fact, 
the example above shows  that the reflexivity of  $\cE \otimes \cF$  does
in general not imply that $\cE$ and  $\cF$ are reflexive too.
\end{remark} 

As preparation for our last statement we need a  result which is analogous
to Lemma \ref{Huneke-U-lifting}. 

\begin{lemma} \label{lifting_for sheaves}
Let $M$ denote a reflexive $R$-module and let $f \in R$ be a homogeneous
$M$-regular element. Suppose that $H^0(M/f M)$ has depth $t \geq 3$. Then
$M$ has depth $t+1$ and $H^0(M/f M) = M/f M$.
\end{lemma}

\begin{proof} Multiplication by $f$ provides the exact cohomology sequence 
\begin{eqnarray*}
\lefteqn{0 \to \HH^0(M/f M) \to \HH^1(M)(-k) \stackrel{f}{\longrightarrow} \HH^1(M) \to
\HH^1(M/f M) \to \hspace*{4cm} } \\ 
& & \hspace*{3cm}  \HH^2(M)(-k) \stackrel{f}{\longrightarrow} \HH^2(M) \to
\HH^2(M/f M) \to  \ldots.
\end{eqnarray*} 
Since we have 
$$
\HH^2(M/f M) \cong \HH^2(H^0(M/f M)) = 0 
$$ 
by assumption on the depth $t$ the multiplication by $f$ is surjective on 
$\HH^2(M)$. But $\HH^2(M)$ is a finitely generated $R$-module because $M$
is reflexive. Hence Nakayama's lemma implies $\HH^2(M) = 0$. Furthermore,
the reflexivity of $M$ provides $\HH^1(M) = 0$.  Thus, considering the
exact sequence above we conclude that 
$$
\HH^0(M/f M) = \HH^1(M/f M) = 0. 
$$ 
It follows that $H^0(M/fM) \cong M/f M$ and we are done. 
\end{proof}
 
For short we say that a coherent sheaf $\cF$ on $\PP^n$ splits if it is a
direct sum of line bundles. This holds true if and only if $H^0_*(\cF)$ is
a finitely generated free $R$-module.  

\begin{proposition} Let $\cE,\cF$ denote reflexive
  sheaves on $\PP^n$ ($n \geq 1$). 
\begin{itemize} 
\item[(a)]  If $\cE \otimes \cF$ splits and $\cE$ and $\cF$
intersect very properly then $\cE$ and $\cF$ split. 
\item[(b)] If $\cE$ is a vector bundle such that  $H^1_*(\cE \otimes \cE^*) = 0$ then
$\cE$ splits.
\end{itemize} 
\end{proposition} 

\begin{proof} (a) The assumption and Proposition
 \ref{proper_makes_tensor_product_nice} provide that the tensor product of 
$E = H^0_*(\cE)$ and $F = H^0_*(\cF)$ is a free $R$-module. Therefore, due
to Theorem \ref{CM-lifting_theorem}, $E$ 
and $F$ are Cohen-Macaulay, thus free. 

(b) Put again $E = H^0_*(\cE)$. Our assumption implies $\depth E \tr E^*
\geq 3$. Since $E \tr E^*) \cong (E \tk E^*)/ \Delta (E \tk E^*$
we conclude by Lemma \ref{lifting_for sheaves} that $\depth E \tk E^* \geq
n+4$ and in particular $\HH^{n+2} (E \tk E^*) = 0$. 

Let us denote the $K$-dual $\Hom_K(M, K)$ of a graded $R$-module $M$ by
$M^{\vee}$. The K\"unneth formulas and local duality provide  
\begin{eqnarray*} 
\HH^{n+2} (E \tk E^*) & = & \bigoplus_{i+j=n+2} \HH^i(E) \tk \HH^j(E^*) \\
& = & \bigoplus_{i=2}^n \HH^i(E) \tk \HH^i(E)^{\vee}(-n-1). 
\end{eqnarray*} 
Since each direct summand must be trivial we conclude that $\HH^i(E) = 0$ 
for all $i \leq n$, i.e., $E$ is a maximal Cohen-Macaulay module and thus
free.  
\end{proof} 

Part (b) of the last statement has been proved over the complex numbers by
Luk and Yau in \cite{Luk_Yau} using analytic tools. The general case has
been first proved by Huneke and Wiegand (\cite{Huneke-Wiegand-Math.Scand},
Theorem 5.2) using their results in \cite{Huneke-Wiegand-Math.Ann}. 
\medskip

An analysis of the proof of claim (b) shows that its conclusion follows once
we know that  $[\HH^{n+2}(E \tk E^*)]_{n+1} = 0$.  Thus it seems
conceivable that the hypothesis could be weakened  by assuming only the
vanishing of $H^1(\cE \otimes \cE^*(t))$ for a {\bf finite} number of
degrees $t$. We leave this as an open problem. 


\end{document}